\documentclass[10pt]{amsart}

\usepackage{enumerate, amsmath, amsfonts, amssymb, amsthm, mathtools, thmtools, wasysym, graphics, graphicx, xcolor, frcursive,xparse,comment,bbm}
\usepackage[all]{xy}
\usepackage[colorlinks=true, pdfstartview=FitV, linkcolor=blue, citecolor=blue, urlcolor=blue]{hyperref}

\usepackage{url, hypcap}
\hypersetup{colorlinks=true, citecolor=darkblue, linkcolor=darkblue}
\usepackage{mathrsfs}

\usepackage{tikz}
\usetikzlibrary{arrows,automata}
\usetikzlibrary{patterns,snakes}
\usetikzlibrary{shapes.geometric,calc}
\tikzstyle{arrow}=[thick, ->, >=stealth]
\tikzset{
	edge/.style={->,> = latex'}
}
\usetikzlibrary{calc,through,backgrounds,shapes,matrix}
\usepackage{graphicx}

\usepackage[letterpaper,margin=1.25in]{geometry}

\definecolor{darkblue}{rgb}{0.0,0,0.7} 
\definecolor{lightblue}{rgb}{0,135,147}

\definecolor{red}{rgb}{1,0,0}

\definecolor{darkred}{rgb}{0.7,0,0} 
\definecolor{lightgrey}{rgb}{0.7,0.7,0.7} 

\newtheorem{theorem}{Theorem}[section]
\newtheorem{proposition}[theorem]{Proposition}

\theoremstyle{definition}
\newtheorem{definition}[theorem]{Definition}
\newtheorem{example}[theorem]{Example}

\newtheorem{remark}[theorem]{Remark}

\numberwithin{equation}{section}

\definecolor{darkred}{rgb}{0.7,0,0} 
\newcommand{\defn}[1]{{\color{darkred}\emph{#1}}} 

\usepackage[colorinlistoftodos]{todonotes}

\title
  {Normal distributions of finite Markov chains}
  
\author[J.~Rhodes]{John Rhodes}
\address[J. Rhodes]{Department of Mathematics, University of California, Berkeley, CA 94720, U.S.A.}
\email{rhodes@math.berkeley.edu, blvdbastille@gmail.com}

\author[A.~Schilling]{Anne Schilling}
\address[A. Schilling]{Department of Mathematics, UC Davis, One Shields Ave., Davis, CA 95616-8633, U.S.A.}
\email{anne@math.ucdavis.edu}

\date{\today}

\begin{document}

\begin{abstract}
We show that the stationary distribution of a finite Markov chain can be expressed as the sum of certain
normal distributions. These normal distributions are associated to planar graphs consisting of a straight line
with attached loops. The loops touch only at one vertex either of the straight line or of another attached loop.
Our analysis is based on our previous work, which derives the stationary distribution of a finite Markov chain
using semaphore codes on the Karnofsky--Rhodes and McCammond expansion of the right Cayley graph
of the finite semigroup underlying the Markov chain. 
\end{abstract}

\maketitle

\section{Introduction}
In our previous paper~\cite{RS.2017}, we developed a general theory to compute the stationary distribution of
a finite Markov chain. Every finite state Markov chain $\mathcal{M}$ has a random letter representation, that is, a 
representation of a semigroup $S$ acting on the left on the state space $\Omega$~\cite{LPW.2009}. Combining the
Karnofsky--Rhodes and the McCammond expansion of the right Cayley graph of $S$, we were able to provide a construction
of the stationary distribution using finite semigroup theory without the use of linear algebra. The construction relies on the 
concept of lumping; the distributions for the expanded graphs can be computed thanks to normal forms of the elements. 
The stationary distribution of the original Markov chain $\mathcal{M}$ is then obtained by lumping.

In this paper, we show that the stationary distribution of any finite Markov chain can be obtained from
certain \defn{normal} (or \defn{Gau\ss ian}) \defn{distributions}. The normal distributions are derived from planar graphs 
by adding directed loops (or circles) to the straight line, which only touch the graph at one point. Let us outline the construction 
of these normal forms in the remainder of the introduction.

\subsection{Straight line}
\label{section.straight}
We start with a straight line starting at $\mathbbm{1}$ with $n$ further vertices:
\begin{center}
\begin{tikzpicture}[auto]
\node (A) at (0, 0) {$\mathbbm{1}$};
\node (B) at (1.5,0) {$1$};
\node(C) at (3,0) {$2$};
\node(Cp) at (4.5,0) {};
\node(Dp) at (6,0) {};
\node(F) at (5.3,0) {$\cdots$};
\node(D) at (7.5,0) {$n-1$};
\node(E) at (9.2,0) {$n$};
\draw[edge,thick] (A) -- (B);
\draw[edge,thick] (B) -- (C);
\draw[edge,thick] (C) -- (Cp);
\draw[edge,thick] (Dp) -- (D);
\draw[edge,thick] (D) -- (E);
\end{tikzpicture}
\end{center}

\subsection{Adding loops}

A \defn{loop} is a sequence of vertices connected by edges $v_0 \longrightarrow v_1 \longrightarrow \cdots \longrightarrow
v_k$ such that $v_0 = v_k$, but all other vertices $v_i$ with $0\leqslant i<k$ are distinct.

Add a loop $\ell$ to any vertex of the straight line constructed in Section~\ref{section.straight} (except $\mathbbm{1}$) 
with $k\geqslant 0$ new vertices, which only touches one existing vertex $v$.
\begin{center}
\begin{tikzpicture}[auto]
\node (A) at (0, 0) {$\mathbbm{1}$};
\node (B) at (1.5,0) {$1$};
\node(C) at (3,0) {$2$};
\node(Cp) at (4.5,0) {$v$};
\node(Dp) at (6,0) {};
\node(F) at (5.3,0) {$\cdots$};
\node(D) at (7.5,0) {$n-1$};
\node(E) at (9.2,0) {$n$};
\node(W) at (5,0.5) {$v_1$};
\node (V) at (5.5,1) {$v_2$};
\node (X) at (4.5,3) {$\cdots$};
\node(Y) at (4,0.5) {$v_k$};
\node (Z) at (3.5,1){$v_{k-1}$};

\draw[edge,thick] (A) -- (B);
\draw[edge,thick] (B) -- (C);
\draw[edge,thick] (C) -- (Cp);
\draw[edge,thick] (Dp) -- (D);
\draw[edge,thick] (D) -- (E);
\path (Cp) edge[->,thick, bend right = 30] (W);
\path (W) edge[->,thick, bend right = 30] (V);
\path (V) edge[->,thick, bend right = 60] (X);
\path (X) edge[->,thick, bend right = 60] (Z);
\path (Z) edge[->,thick, bend right = 30] (Y);
\path (Y) edge[->,thick, bend right = 30] (Cp);
\end{tikzpicture}
\end{center}

The cut of $\ell$ is
\begin{center}
\begin{tikzpicture}[auto]
\node (A) at (0, 0) {$v$};
\node (B) at (1.5,0) {$v_1$};
\node(C) at (3,0) {$v_2$};
\node(Cp) at (4.5,0) {};
\node(Dp) at (6,0) {};
\node(F) at (5.3,0) {$\cdots$};
\node(D) at (7.5,0) {$v_{k-1}$};
\node(E) at (9.2,0) {$v_k$};
\draw[edge,thick] (A) -- (B);
\draw[edge,thick] (B) -- (C);
\draw[edge,thick] (C) -- (Cp);
\draw[edge,thick] (Dp) -- (D);
\draw[edge,thick] (D) -- (E);
\end{tikzpicture}
\end{center}

Continue to add loops at any vertex (except $\mathbbm{1}$), including the new vertices. Multiple loops
at a given vertex are allowed.
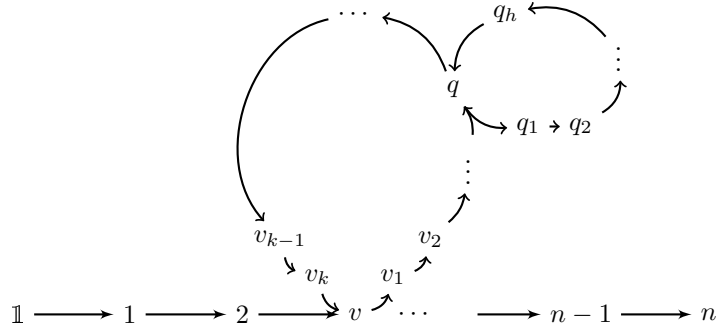
\begin{center}
\begin{tikzpicture}[auto]
\node (A) at (0, 0) {$\mathbbm{1}$};
\node (B) at (1.5,0) {$1$};
\node(C) at (3,0) {$2$};
\node(Cp) at (4.5,0) {$v$};
\node(Dp) at (6,0) {};
\node(F) at (5.3,0) {$\cdots$};
\node(D) at (7.5,0) {$n-1$};
\node(E) at (9.2,0) {$n$};
\node(W) at (5,0.5) {$v_1$};
\node (V) at (5.5,1) {$v_2$};
\node(Vp) at (6,2) {$\vdots$};
\node(Vpp) at (5.8,3) {$q$};
\node (X) at (4.5,4) {$\cdots$};
\node(Y) at (4,0.5) {$v_k$};
\node (Z) at (3.5,1){$v_{k-1}$};
\node(L1) at (6.8,2.5) {$q_1$};
\node(L2) at (7.5,2.5) {$q_2$};
\node(L3) at (8,3.5) {$\vdots$};
\node(L4) at (6.5,4) {$q_h$};

\draw[edge,thick] (A) -- (B);
\draw[edge,thick] (B) -- (C);
\draw[edge,thick] (C) -- (Cp);
\draw[edge,thick] (Dp) -- (D);
\draw[edge,thick] (D) -- (E);
\path (Cp) edge[->,thick, bend right = 30] (W);
\path (W) edge[->,thick, bend right = 30] (V);
\path (V) edge[->,thick, bend right = 20] (Vp);
\path (Vp) edge[->,thick, bend right = 20] (Vpp);
\path (Vpp) edge[->,thick, bend right = 30] (X);
\path (X) edge[->,thick, bend right = 60] (Z);
\path (Z) edge[->,thick, bend right = 30] (Y);
\path (Y) edge[->,thick, bend right = 30] (Cp);
\path (Vpp) edge[->,thick, bend right = 30] (L1);
\path (L1) edge[->,thick] (L2);
\path (L2) edge[->,thick, bend right = 30] (L3);
\path (L3) edge[->,thick, bend right = 30] (L4);
\path (L4) edge[->,thick, bend right = 30] (Vpp);
\end{tikzpicture}
\end{center}
Let $\overline{G}$ be the directed graph obtained by this procedure. Notice that each such $\overline{G}$ can be drawn in 
the plane.

\subsection{Kleene expressions}
\label{section.kleene}
Given a finite alphabet $A$, assign a letter $a\in A$ to each arrow in the graph $\overline{G}$.
The result is called a \defn{loop graph}, denoted $G$.

\begin{example}
\label{example.G}
For the alphabet $A=\{a,b,c,d,x\}$, we might obtain
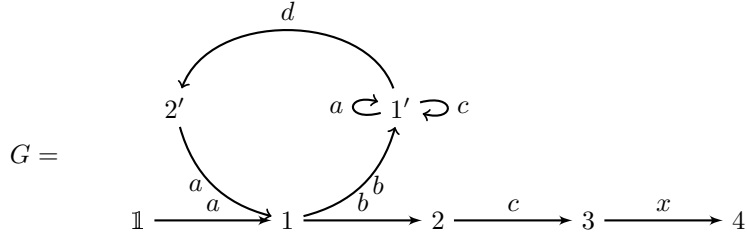
\begin{center}
\raisebox{1cm}{$G =$} \qquad  
\begin{tikzpicture}[auto]
\node (A) at (0, 0) {$\mathbbm{1}$};
\node (B) at (2,0) {$1$};
\node(C) at (4,0) {$2$};
\node(D) at (6,0) {$3$};
\node(E) at (8,0) {$4$};
\node(W) at (3.5,1.5) {$1'$};
\node (V) at (0.5,1.5) {$2'$};

\draw[edge,thick] (A) -- (B) node[midway, above] {$a$};
\draw[edge,thick] (B) -- (C) node[midway, above] {$b$};
\draw[edge,thick] (C) -- (D) node[midway, above] {$c$};
\draw[edge,thick] (D) -- (E) node[midway, above] {$x$};
\path (B) edge[->,thick, bend right = 30] node[midway,right] {$b$} (W);
\path (V) edge[->,thick, bend right = 30] node[midway,left] {$a$} (B);
\path (W) edge[->,thick, bend right = 70] node[midway,above] {$d$} (V);
\path (W) edge [thick,loop right] node {$c$} (W);
\path (W) edge [thick,loop left] node {$a$} (W);
\end{tikzpicture}
\end{center}
\end{example}
In general, this procedure gives a non-deterministic automata since different edges emitting from a vertex can be labeled 
by the same letter. In the above example, vertex 1 has two arrows labeled $b$ coming out of it. 

Denote the set of all paths in a loop graph $G$ starting at $\mathbbm{1}$ and ending at $n$ (the last vertex 
on the initial straight line underlying $G$) by $\mathcal{P}_G$. Here a path is given by
\[
	\mathbbm{1} \stackrel{a_1}{\longrightarrow} v_1 \stackrel{a_2}{\longrightarrow} \cdots \stackrel{a_k}
	{\longrightarrow} v_k=n,
\]
where $v_i$ are vertices in $G$ and $a_i \in A$ are the labels on the edges.

There is a simple inductive way to describe $\mathcal{P}_G$ using \defn{Kleene expressions}.
Given a set $L$, define $L^0 = \{\varepsilon \}$ given by the empty string, $L^1 = L$, and recursively
$L^{i+1} = \{wa \mid w \in L^i, a \in L\}$ for each integer $i>0$. Then the \defn{Kleene star} is
\[
	L^\star = \bigcup_{i\geqslant 0} L^i.
\]
A Kleene expression only involves letters in $A$, unions, and $\star$.
To obtain a Kleene expression for $\mathcal{P}_G$, perform the following doubly recursive procedure:

\smallskip
\noindent
\textbf{Algorithm 1.}

\smallskip
\noindent
\textbf{Induction basis:} Start at vertex $\mathbbm{1}$ and with the empty expression $L$.

\smallskip
\noindent
\textbf{Induction step:} Suppose one is at vertex $i\neq n$ (or $\mathbbm{1}$) on the straight line path underlying $G$.
\begin{enumerate}
\item Continue to the next vertex $i+1$ (or $1$) on the straight line path underlying $G$ and append the label $a$
on the edge from $i \stackrel{a}{\longrightarrow} i+1$ (or $\mathbbm{1} \stackrel{a}{\longrightarrow} 1$) to $L$.
\item If there are loops $\ell_1,\ell_2,\ldots,\ell_k$ at vertex $i+1$ (or $1$), append the formal expression
\[
	\{\ell_1,\ell_2,\ldots,\ell_k\}^\star
\]
to $L$. The loops $\ell_1,\ell_2,\ldots,\ell_k$ are in one-to-one correspondence with the edges coming into vertex
$i+1$.
\item If $i+1\neq n$, continue with the next induction step. Else stop and output $L$.
\end{enumerate}

\smallskip
\noindent
\textbf{Algorithm 2.} For each symbol $\ell_i$ in the expression for $L$, do the following:
\begin{enumerate}
\item Consider the loop $\ell_i = \left( v_0 \stackrel{a_1}{\longrightarrow} v_1 \stackrel{a_2}{\longrightarrow} \cdots 
\stackrel{a_k}{\longrightarrow} v_k=v_0 \right)$ from vertex $v_0$ to $v_0$ in $G$. Consider the subgraph
of $G$ with straight line $v_1 \stackrel{a_2}{\longrightarrow} \cdots \stackrel{a_k}{\longrightarrow} v_k$ and all
further loops that are attached to any of the vertices $v_i$ in $G$. Attach $\mathbbm{1}$ to $v_1$. The resulting
graph $G^{(i)}$ is a new loop graph. Perform Algorithm 1 on $G^{(i)}$ to obtain a Kleene expression $L^{(i)}$.
Replace the symbol $\ell_i$ in $L$ by $L^{(i)}$.
\item Continue this process until $L$ does not contain any further expressions $\ell_i$ for some loop $\ell_i$,  that
is, $L$ only contains unions, $\star$ and elements in the alphabet $A$.
Then the Kleene expression for $\mathcal{P}_G$ is $L$.
\end{enumerate}
The resulting expressions can be made into unionless expressions by using \defn{Zimin words}
\begin{equation}
\label{equation.zimin}
	\{a\}^\star = a^\star \qquad \text{and} \qquad \{a, b\}^\star = (a^\star b)^\star a^\star \qquad \text{for $a,b\in A$.}
\end{equation}
Expressions for larger unions can be obtained by induction using~\eqref{equation.zimin}.

\begin{example}
Let $G$ be as in Example~\ref{example.G}. Then
\[
	L = a \ell_1^\star bcx,
\]
where $\ell_1$ is the loop attached to vertex 1. Cut this loop and continue the process to obtain
\[
	\ell_1 = b \{\ell_1',\ell_2'\}^\star da,
\]
where $\ell_1'$ is the loop at vertex $1'$ labelled $a$ and $\ell_2'$ is the loop at vertex $1'$ labelled $c$.
We have $\ell_1'=a$ and $\ell_2'=c$, so that altogether we find
\[
	L = a (b\{a,c\}^\star d a)^\star bcx = a (b(a^\star c)^\star a^\star d a)^\star bcx,
\]
where in the last step we used the Zimin words to get rid of the unions. This is a Kleene expression for 
$\mathcal{P}_G$.
\end{example}
See Example~\ref{example.pict1} for another example and also compare this construction to the definition of 
$\mathsf{Pict}$ in Definition~\ref{definition.pict}.

\subsection*{Main results}
We are now going to define normal distributions. 

\begin{definition}[Normal distribution]
\label{definition.normal}
Let $G$ be a loop graph with edges labeled by letters in the alphabet $A$. Associate the indeterminate $x_a$ to
$a\in A$. Then the \defn{normal distribution} of $G$ is defined as
\[
	\Psi_G = \sum_{p \in \mathcal{P}_G} \prod_{a\in p} x_a.
\]
\end{definition}

We may use the Kleene expressions of the previous section for $\mathcal{P}_G$.
The advantage in doing so is that one can immediately obtain rational expressions. Namely, using the geometric series, 
we find that
\[
	\sum_{s\in a^\star} \prod_{i \in s} x_i = \sum_{\ell=0}^\infty x_a^\ell = \frac{1}{1-x_a}.
\]
Similarly
\[
	\sum_{s \in \{a,b\}^\star} \prod_{i\in s} x_i = \sum_{s \in a^\star (ba^\star)^\star} \prod_{i\in s} x_i
	= \frac{1}{1-x_a} \cdot \frac{1}{1-\frac{x_b}{1-x_a}}
	= \frac{1}{1-x_a-x_b}.
\]
In general, using the recursion~\eqref{equation.zimin} we derive by induction
\begin{equation}
\label{equation.geometric}
	\sum_{s \in \{a_1,a_2,\ldots,a_n\}^\star} \prod_{i\in s} x_i = \frac{1}{1-x_{a_1} - x_{a_2} - \cdots - x_{a_n}}.
\end{equation}

Our main theorem is the following. 
\begin{theorem}
\label{theorem.main}
The stationary distribution $\Psi^{\mathcal{M}}$ of a finite Markov chain $\mathcal{M}$ is the sum of normal 
distributions $\Psi_G$ or certain limits of $\Psi_G$, where $G$ is a loop graph.
\end{theorem}

The proof of Theorem~\ref{theorem.main} is given in Section~\ref{section.proof}. A more precise version
of Theorem~\ref{theorem.main} is stated in Theorem~\ref{theorem.main 1}.

The paper is outlined as follows. In Section~\ref{section.walks}, we review the main results from~\cite{RS.2017},
in particular the expressions for the stationary distribution of a finite Markov chain in terms of semaphore codes
of the Karnofsky--Rhodes expansion of the right Cayley graph of the underlying semigroup. In Section~\ref{section.normal},
we review the McCammond expansion and its relation to semaphore codes and provide the definition of $\mathsf{Pict}$.
The map $\mathsf{Pict}$ is used to give a proof of Theorem~\ref{theorem.main}.
The original definition of $\mathsf{Pict}$ is due to McCammond, but the applications to random walks are due to the 
authors.

\subsection*{Acknowledgments}
We are grateful to Jon McCammond and Ben Steinberg for discussions.
The map $\mathsf{Pict}$ of Definition~\ref{definition.pict} is due to McCammond, told to the first author in 1994, 
written by the first author in 2008, and simplified here.

The first author thanks the Simons Foundation Collaboration Grants for Mathematicians for travel grant \#313548.
The second author was partially supported by NSF grants DMS--1760329 and DMS--1764153. 

\section{Stationary distributions of Markov chains}
\label{section.walks}

In this section, we provide definitions and review the necessary results we need from~\cite{RS.2017}.

\subsection{Markov chains}
\label{section.markov}

A \defn{Markov chain} $\mathcal{M}$ consists of a finite or countable state space $\Omega$ together with transition
probabilities $\mathcal{T}_{s',s}$ for the transition $s \longrightarrow s'$ for $s,s' \in \Omega$.
The matrix $\mathcal{T} = (\mathcal{T}_{s',s})_{s,s'\in \Omega}$ is called the \defn{transition matrix},
which is a column-stochastic matrix, meaning that the column sums of $\mathcal{T}$ are equal to one.

A Markov chain is \defn{irreducible} if for any $s,s' \in \Omega$ there exists an integer $m$ (possibly depending
on $s$, $s'$) such that $\mathcal{T}_{s',s}^m >0$. In other words, one can get from any state $s$ to any other state $s'$ 
using only steps with positive probability. 
A state $s\in \Omega$ is called \defn{recurrent} if the system returns
to $s$ in finitely many steps with probability one.

The \defn{stationary distribution} of $\mathcal{M}$ is a vector $\Psi = (\Psi_s)_{s\in \Omega}$ such that 
$\mathcal{T} \Psi = \Psi$ and $\sum_{s \in \Omega} \Psi_s=1$. In other words, $\Psi$ is a right-eigenvector of 
$\mathcal{T}$ with eigenvalue one. If the Markov chain is irreducible, the stationary distribution is unique~\cite{LPW.2009}.

Next we define \defn{lumping} of Markov chains. Partition the state space $\Omega$
into $(\Omega_1,\ldots,\Omega_\ell)$ such that
\[
	\Omega_i \cap \Omega_j = \emptyset \quad \text{for $i\neq j$ and} \quad
	\Omega = \bigcup_{i=1}^\ell \Omega_i.
\]
One may view such a partition as an equivalence relation $s\sim s'$ if $s,s'\in \Omega_i$ for some $1\leqslant i \leqslant \ell$.
We say that $\mathcal{M}$ can be lumped with respect to the partition $(\Omega_1,\ldots,
\Omega_\ell)$ if the transition matrix $\mathcal{T}$ satisfies~\cite[Lemma~2.5]{LPW.2009} \cite{KS.1976}
for all $1\leqslant i,j \leqslant \ell$
\begin{equation}
\label{equation.lumping}
	\sum_{t\in \Omega_j} \mathcal{T}_{t,s} = \sum_{t\in \Omega_j} \mathcal{T}_{t,s'}
	\qquad \text{for all $s,s' \in \Omega_i.$}
\end{equation}
The lumped Markov chain is a random walk on the equivalence classes, whose stationary distribution labeled by
$w$ is $\sum_{s \sim w} \Psi_s$.

Every finite state Markov chain $\mathcal{M}$ has a random letter representation, that is, a representation of a 
semigroup $S$ acting on the left on the state space $\Omega$ (see~\cite[Proposition 1.5]{LPW.2009} 
and~\cite[Theorem 2.3]{ASST.2015}). In this setting, we transition $s \stackrel{a}{\longrightarrow} s'$ with probability 
$0\leqslant x_a\leqslant 1$, where $s, s'\in \Omega$, $a\in S$ and $s'=a.s$ is the action of $a$ on the state $s$. 
Let $A=\{a\in S \mid x_a>0\}$. We assume that $A$ generates $S$; if not, it suffices to consider the subsemigroup
generated by $A$. Note that $\sum_{a\in A} x_a =1$. The \defn{transition matrix} $\mathcal{T}$ of $\mathcal{M}$ is the 
$|\Omega| \times |\Omega|$-matrix 
\begin{equation}
\label{equation.transition matrix}
	\mathcal{T}_{s',s} = \sum_{\substack{a \in A\\ s \stackrel{a}{\longrightarrow} s'}} x_a
	\qquad \text{for $s,s' \in \Omega$.}
\end{equation}
Note that we may assume that the action of $S$ on $\Omega$ is faithful as this does not affect the random
walk.

If $S$ is a semigroup, then $S^{\mathbbm{1}}$ denotes $S$ with 
an adjoint identity $\mathbbm{1}$ even if $S$ already has an identity.

\begin{definition}[Ideal]
Let $S$ be a semigroup. A two-sided \defn{ideal} $I$ (or ideal for short) is a subset $I \subseteq S$ such that 
$u I v \subseteq I$ for all $u,v \in S^{\mathbbm{1}}$. Similarly, a \defn{left ideal} $I$ is a subset 
$I \subseteq S^{\mathbbm{1}}$ such that $u I \subseteq I$ for all $u\in S^{\mathbbm{1}}$. 
\end{definition}

If $I,J$ are ideals of $S$, then $IJ \subseteq I \cap J$, so that $I \cap J \neq \emptyset$. Hence every
finite semigroup has a unique minimal ideal denoted $K(S)$.
As shown in~\cite{Clifford.Preston.1961,KRT.1968}, the minimal ideal $K(S)$ of a finite semigroup $S$ is the disjoint 
union of all the minimal left ideals of $S$ and the Rees Theorem applies. By~\cite[Remark 2.8]{ASST.2015} 
the faithful left action of $S$ on $\Omega$ is isomorphic to the left action of $S$ on $K(S)$.

Let $(S,A)$ be a semigroup $S$ together with a choice of generators $A$ for $S$.
Define \defn{$\mathcal{M}(S,A)$} to be the Markov chain, where the transition
$s \stackrel{a}{\longrightarrow} s'$ for $s,s'\in S$ and $a\in A$ is given by $s'=as$ in the left Cayley graph
with probability $0<x_a\leqslant 1$. Note that we are assuming that all probabilities $x_a$ for $a\in A$
are nonzero. Then it was shown in~\cite{HM.2011} (see also~\cite[Proposition 3.2]{ASST.2015}) that the
recurrent states of $\mathcal{M}(S,A)$ are the elements in $K(S)$. Furthermore, the connected components
of the recurrent states in the random walk are the minimal left ideals of $S$. The restriction of the random walk 
to any minimal left ideal is irreducible. Moreover, the chain so obtained is independent of the chosen minimal left ideal.
This random walk and the random walk with states a left ideal $L$ of $K(S)$ and $S$ acting on the left made faithful, 
that is $x \stackrel{a}{\longrightarrow} y$ for $x \in L$ and $y = ax$, are essentially the same. So we may not distinguish 
the two cases. 

\subsection{Karnofsky--Rhodes expansion}
\label{section.KR}

In this section, we define the right Cayley graph of a finite semigroup and its Karnofsky--Rhodes expansions.

\begin{definition}[Right Cayley graph]
Let $(S,A)$ be a finite semigroup $S$ together with a set of generators $A$. 
The \defn{right Cayley graph} $\mathsf{RCay}(S,A)$ of $S$ with respect to $A$ is the rooted graph with
vertex set $S^{\mathbbm{1}}$, root $r=\mathbbm{1} \in S^{\mathbbm{1}}$, and edges 
$s \stackrel{a}{\longrightarrow} s'$ for all $(s,a,s') \in S^{\mathbbm{1}} \times A \times S^{\mathbbm{1}}$, where $s'=sa$ 
in $S^{\mathbbm{1}}$.
\end{definition}

A \defn{path} $p$ in $\mathsf{RCay}(S,A)$ is a sequence 
\[
	p = \left(v_1 \stackrel{a_1}{\longrightarrow} \cdots \stackrel{a_\ell} {\longrightarrow} v_{\ell+1} \right),
\]
where $v_i\in S^{\mathbbm{1}}$ are vertices in $\mathsf{RCay}(S,A)$ and $v_i \stackrel{a_i}{\longrightarrow} v_{i+1}$ are
edges in $\mathsf{RCay}(S,A)$. The endpoint of $p$ is $\tau(p) := v_{\ell+1}$. The length of the path $p$
is $\ell(p):=\ell$, which equals the number of edges. 
A \defn{simple path} is a path that does not visit any vertex twice. Empty paths are considered simple.
A path which starts and ends at the same vertex is called a \defn{circuit}. A circuit that is simple, when the last vertex is 
removed, is called a \defn{loop}.

\begin{definition}[Transition edges]
An edge $s \stackrel{a}{\longrightarrow} s'$ in the right Cayley graph $\mathsf{RCay}(S,A)$
is a \defn{transition edge} if there is no directed path from $s'$ to $s$ in $\mathsf{RCay}(S,A)$.
In other words, there does not exist any sequence $a_1,\ldots,a_k \in A$ with $k\geqslant 1$ such that 
$s' (a_1 \cdots a_k) = s$.
\end{definition}

Let us now define the Karnofsky--Rhodes expansion of the right Cayley graph (see
also~\cite[Definition~4.15]{MRS.2011} and~\cite[Section 3.4]{MSS.2015}).
Let $(A^+,A)$ be the free semigroup with generators $A$, where $A^+$ is the set of all words $a_1 \ldots a_\ell$ of 
length $\ell \geqslant 1$ over $A$ with multiplication given by concatenation. When we write 
$[a_1 \cdots a_\ell]_S$, we mean the element in $S$ when taking the product in the semigroup of the generators
$a_i\in A$.

\begin{definition}[Karnofksy--Rhodes expansion]
The \defn{Karnofsky--Rhodes expansion} $\mathsf{KR}(S,A)$ is obtained as follows. Start with the right Cayley graph 
$\mathsf{RCay}(A^+,A)$. Identify two paths in $\mathsf{RCay}(A^+,A)$ 
\begin{equation*}
	p := \left( \mathbbm{1} \stackrel{a_1}{\longrightarrow} v_1 \stackrel{a_2}{\longrightarrow} \cdots \stackrel{a_\ell}
	{\longrightarrow} v_\ell \right)
	\quad \text{and} \quad 
	p' := \left( \mathbbm{1} \stackrel{a'_1}{\longrightarrow} v'_1 \stackrel{a'_2}{\longrightarrow} \cdots 
	\stackrel{a'_{\ell'}}{\longrightarrow} v'_{\ell'} \right)
\end{equation*}
in $\mathsf{KR}(S,A)$ if and only if the corresponding paths in $\mathsf{RCay}(S,A)$
\begin{equation*}
	[p]_S := \left( \mathbbm{1} \stackrel{a_1}{\longrightarrow} [v_1]_S \stackrel{a_2}{\longrightarrow} \cdots 
	\stackrel{a_\ell}{\longrightarrow} [v_\ell]_S \right)
	\quad \text{and} \quad 
	[p']_S := \left( \mathbbm{1} \stackrel{a'_1}{\longrightarrow} [v'_1]_S \stackrel{a'_2}{\longrightarrow} \cdots 
	\stackrel{a'_{\ell'}}{\longrightarrow} [v'_{\ell'}]_S \right),
\end{equation*}
where $v_i=a_1 a_2 \ldots a_i$ and $v_i' = a_1' a_2' \ldots a'_i$, end at the same vertex $[v_\ell]_S = [v'_{\ell'}]_S$ 
and in addition the set of transition edges of $[p]_S$ and $[p']_S$ in $\mathsf{RCay}(S,A)$ is equal. 
\end{definition}

\begin{example}
\label{example.KR}
Consider the right Cayley graph of the Klein $4$-group $Z_2 \times Z_2$ with zero with generators 
$\{a,b,\square\}$, where  $a=(1,-1)$, $b=(-1,1)$, and $\square$ is the zero. The right Cayley graph
$\mathsf{RCay}(Z_2 \times Z_2 \cup \{\square\},\{a,b,\square\})$ is
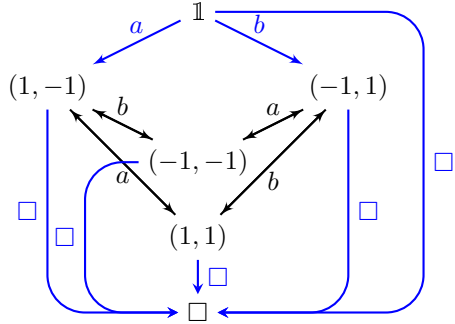
\begin{center}
\begin{tikzpicture}[auto]
\node (A) at (0, 0) {$\mathbbm{1}$};
\node (B) at (-2,-1) {$(1,-1)$};
\node(C) at (2,-1) {$(-1,1)$};
\node(D) at (0,-2) {$(-1,-1)$};
\node(E) at (0,-3) {$(1,1)$};
\node(F) at (0,-4) {$\square$};
\draw[edge,blue,thick] (A) -- (B) node[midway, above] {$a$};
\draw[edge,blue,thick] (A) -- (C) node[midway, above] {$b$};
\draw[edge,thick] (B) -- (D) node[midway, above] {$b$};
\draw[edge,thick] (D) -- (B);
\draw[edge,thick] (C) -- (D) node[midway, above] {$a$};
\draw[edge,thick] (D) -- (C);
\draw[edge,thick] (B) -- (E) node[midway, below] {$a$};
\draw[edge,thick] (E) -- (B);
\draw[edge,thick] (C) -- (E) node[midway, below] {$b$};
\draw[edge,thick] (E) -- (C);
\draw[arrow, rounded corners=5mm,blue] (A) -- ($(A) + (3,0)$) -- ($(C)+(1,0)$)  node[pos=2, right]{$\square$}
-- ($(F) + (3,0)$)  -- (F); 
\draw[arrow, rounded corners=5mm,blue] (C) -- ($(F) + (2,0)$) node[pos=0.5, right]{$\square$} -- (F); 
\draw[arrow, rounded corners=5mm,blue] (B) -- ($(F) + (-2,0)$) node[pos=0.5, left]{$\square$} -- (F); 
\draw[arrow, rounded corners=5mm,blue] (D) -- ($(D) + (-1.5,0)$) -- ($(F) +(-1.5,0)$) node[pos=0.5, left]{$\square$} -- (F) ; 
\draw[arrow,blue] (E) -- (F) node[midway,right]{$\square$}; 
\end{tikzpicture}
\end{center}
where all three arrows $a,b,\square$ fix the vertex $\square$ at the bottom.
Transition edges are indicated in blue. Double edges mean that right multiplication by the label for either vertex yields 
the other vertex. The Karnofsky--Rhodes expansion of this right Cayley graph is given by
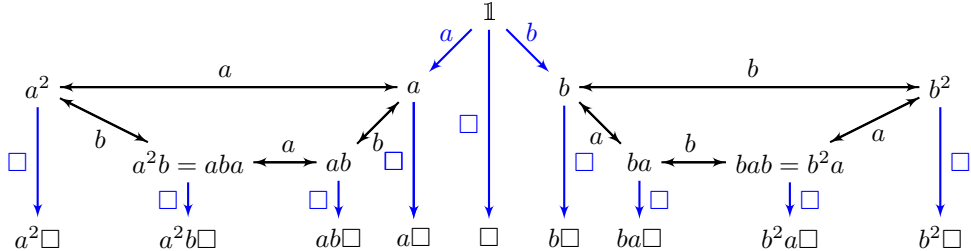
\begin{center}
\begin{tikzpicture}[auto]
\node (A) at (0, 0) {$\mathbbm{1}$};
\node (B) at (-1,-1) {$a$};
\node(C) at (1,-1) {$b$};
\node(D) at (-2,-2) {$ab$};
\node(E) at (2,-2) {$ba$};
\node(F) at (-6,-1) {$a^2$};
\node(G) at (6,-1) {$b^2$};
\node(H) at (-4,-2) {$a^2b=aba$};
\node(I) at (4,-2) {$bab=b^2a$};
\node(ZA) at (0,-3) {$\square$};
\node(ZB) at (-1,-3) {$a\square$};
\node(ZD) at (-2,-3) {$ab\square$};
\node(ZH) at (-4,-3) {$a^2b\square$};
\node(ZF) at (-6,-3) {$a^2\square$};
\node(ZC) at (1,-3) {$b\square$};
\node(ZE) at (2,-3) {$ba\square$};
\node(ZI) at (4,-3) {$b^2a\square$};
\node(ZG) at (6,-3) {$b^2\square$};
\draw[edge,blue,thick] (A) -- (B) node[midway, above] {$a$\;};
\draw[edge,thick,blue] (A) -- (C) node[midway, above] {\;$b$};
\draw[edge,thick] (B) -- (F) node[midway,above] {$a$\;};
\draw[edge,thick] (F) -- (B);
\draw[edge,thick] (B) -- (D) node[midway,below] {$b$};
\draw[edge,thick] (D) -- (B);
\draw[edge,thick] (F) -- (H) node[midway,below] {$b$\;};
\draw[edge,thick] (H) -- (F);
\draw[edge,thick] (D) -- (H) node[midway,above] {\;$a$};
\draw[edge,thick] (H) -- (D);
\draw[edge,thick] (C) -- (E) node[midway,below] {$a$\;};
\draw[edge,thick] (E) -- (C);
\draw[edge,thick] (E) -- (I) node[midway,above] {$b$\;};
\draw[edge,thick] (I) -- (E);
\draw[edge,thick] (C) -- (G) node[midway,above] {\;$b$};
\draw[edge,thick] (G) -- (C);
\draw[edge,thick] (G) -- (I) node[midway,below] {\;$a$};
\draw[edge,thick] (I) -- (G);
\draw[edge,blue,thick] (A) -- (ZA) node[midway, left] {$\square$};
\draw[edge,blue,thick] (B) -- (ZB) node[midway, left] {$\square$};
\draw[edge,blue,thick] (D) -- (ZD) node[midway, left] {$\square$};
\draw[edge,blue,thick] (H) -- (ZH) node[midway, left] {$\square$};
\draw[edge,blue,thick] (B) -- (ZB) node[midway, left] {$\square$};
\draw[edge,blue,thick] (F) -- (ZF) node[midway, left] {$\square$};
\draw[edge,blue,thick] (C) -- (ZC) node[midway, right] {$\square$};
\draw[edge,blue,thick] (E) -- (ZE) node[midway, right] {$\square$};
\draw[edge,blue,thick] (I) -- (ZI) node[midway, right] {$\square$};
\draw[edge,blue,thick] (G) -- (ZG) node[midway, right] {$\square$};
\end{tikzpicture}
\end{center}
where arrows $a,b,\square$ fix all the vertices at the bottom.
\end{example}

\begin{proposition}\cite[Proposition 2.15]{RS.2017}
\label{proposition.KR Cayley}
$\mathsf{KR}(S,A)$ is the right Cayley graph of a semigroup, also denoted by $\mathsf{KR}(S,A)$.
\end{proposition}

\subsection{Stationary distribution}
\label{section.stationary}

We now review the main results of~\cite{RS.2017}, which give the stationary distribution for any
Markov chain $\mathcal{M}(S,A)$ for a finite semigroup with chosen generators $(S,A)$. 
Recall that $\mathcal{M}(S,A)$ is the random walk on the unique minimal ideal $K(S)$ of $S$. More precisely,
the random walk is given by the left action of $S$ on $K(S)$.

To state our results for the stationary distribution, we first need to review the \defn{semaphore codes} associated to 
$(S,A)$~\cite{BPR.2010}. The semaphore code $\mathcal{S}(S,A)$ is the set of all words $a_1 a_2  \ldots a_\ell \in A^+$ 
such that $[a_1 a_2 \cdots a_\ell]_S \in K(S)$, but  $[a_1 a_2 \cdots a_{\ell-1}]_S \not \in K(S)$.

The main results are the following.

\begin{theorem} \cite[Corollary 2.28]{RS.2017}
\label{theorem.stationary real}
The Markov chain $\mathcal{M}(S,A)$ is the lumping of $\mathcal{M}(\mathsf{KR}(S,A))$ with 
stationary distribution
\[
	\Psi_w^{\mathcal{M}(S,A)} = \sum_{\substack{v \in \mathsf{KR}(S,A)\\ [v]_S=w}} \; 
	\Psi_v^{\mathcal{M}(\mathsf{KR}(S,A))} \qquad \text{for all $w \in (S,A)$.}
\]
\end{theorem}

The next result is non-trivial. It requires the assumption that the minimal ideal $K(S)$ is \defn{left zero}, 
that is, $xy = x$ for all $x,y \in K(S)$.

\begin{theorem} \cite[Theorem 2.12]{RS.2017}
\label{theorem.stationary general}
If $K(S)$ is left zero, the stationary distribution of the Markov chain $\mathcal{M}(\mathsf{KR}(S,A))$ is given by
\begin{equation*}
	\Psi^{\mathcal{M}(\mathsf{KR}(S,A))}_w 
	= \sum_{\substack{s \in \mathcal{S}(S,A)\\ [s]_{\mathsf{KR}(S,A)}=w}}\; \prod_{a\in s} x_a
	\qquad \text{for all $w \in K(\mathsf{KR}(S,A))$.}
\end{equation*}
\end{theorem}

As outlined in~\cite[Section 2.9]{RS.2017}, the case when $K(S)$ is not left zero can be constructed
from the case when $K(S)$ is left zero using the flat operation. That is, one adds an additional generator
$\square$ to the alphabet $A$, which acts as zero. The associated probability is $x_\square$. The elements in the
minimal ideal $K(\mathsf{KR}(S \cup \{\square\}, A \cup \{\square\}))$ are of the form $w\square$, where $w \in
\mathsf{KR}(S,A)$. Since $\square v = \square$ for all $v \in \mathsf{KR}(S,A)$, we indeed have that 
$K(\mathsf{KR}(S \cup \{\square\}, A \cup \{\square\}))$ is left zero and hence Theorem~\ref{theorem.stationary general}
applies. Then~\cite[Corollary~2.33]{RS.2017}
\begin{equation}
\label{equation.limit}
	\Psi^{\mathcal{M}(\mathsf{KR}(S,A))}_w  = \lim_{x_\square \to 0}
	\Psi^{\mathcal{M}(\mathsf{KR}(S \cup \{\square\},A \cup \{\square\}))}_w.
\end{equation}

\section{Normal distributions for random walks}
\label{section.normal}

In this section, we prove Theorem~\ref{theorem.main}. By Theorems~\ref{theorem.stationary real} 
and~\ref{theorem.stationary general} and Equation~\eqref{equation.limit}, the stationary distribution 
$\Psi_w^{\mathcal{M}(S,A)}$ is the sum of terms of the form $\prod_{a\in s} x_a$, where $s\in \mathcal{S}(S,A)$
(or limits of such expressions).
In Section~\ref{section.semaphore}, we will explain how the semaphore code $\mathcal{S}(S,A)$ is related to the 
McCammond expansion $\mathsf{Mc} \circ \mathsf{KR}(S,A)$. In Section~\ref{section.pict}, we will
then define the map $\mathsf{Pict}$ on $\mathsf{Mc} \circ \mathsf{KR}(S,A)$ to deduce that
$\Psi_w^{\mathcal{M}(S,A)}$ is a sum of normal forms. A proof of Theorem~\ref{theorem.main} is
given in Section~\ref{section.proof}. Theorem~\ref{theorem.main 1} is a more precise version of
Theorem~\ref{theorem.main}.

\subsection{The McCammond expansion and semaphore codes}
\label{section.semaphore}

Let us now turn to the McCammond expansion~\cite{McCammond.2001,MRS.2011} of the Karnofsky--Rhodes expansion 
of the right Cayley graph of $(S,A)$. Recall that a \defn{simple path} in $\mathsf{KR}(S,A)$ is a path that does not visit 
any vertex twice. Empty paths are considered simple.

\begin{definition}[McCammond expansion]
\label{definition.mccammond}
The \defn{McCammond expansion} $\mathsf{Mc} \circ \mathsf{KR}(S,A)$ of $\mathsf{KR}(S,A)$ is the graph with vertex set 
$V$, which is the set of simple paths in $\mathsf{KR}(S,A)$. The edges are given by 
\begin{equation*}
\begin{split}
	E := \{(p,a,q) \in V \times A \times V \mid& \quad \tau(q) = \tau(p) a, \;\ell(q) \leqslant \ell(p)+1,\\
	& \quad \text{$q$ is an initial segment of $p$ if $\ell(q) \leqslant \ell(p)$} \}.
\end{split}
\end{equation*}
In other words, if the path $pa$ in $\mathsf{KR}(S,A)$ is simple, then $q = pa$. Otherwise $\tau(pa) = v$ is a vertex of $p$
and then $q$ is the initial segment of $p$ up to and including $v$. 
\end{definition}

\begin{remark}
\label{remark.spanning tree}
Note that $\mathsf{Mc} \circ \mathsf{KR}(S,A)$ has a spanning tree $\mathsf{T}$ with the same vertex set
as $\mathsf{Mc} \circ \mathsf{KR}(S,A)$, but only those edges $(p,a,q) \in E$ such that $\ell(q)=\ell(p)+1$.
\end{remark}

\begin{figure}[t]
\begin{center}
\scalebox{0.9}{
\begin{tikzpicture}[auto]
\node (A) at (0, 0) {$\mathbbm{1}$};
\node (B) at (-4,-1) {$a$};
\node(C) at (4,-1) {$b$};
\node(D) at (-3,-2) {$ab$};
\node(E) at (3,-2) {$ba$};
\node(F) at (-5,-2) {$a^2$};
\node(FF) at (-6,-3) {$a^2b$};
\node(FFF) at (-7,-4) {$a^2ba$};
\node(G) at (5,-2) {$b^2$};
\node(H) at (-2,-3) {$aba$};
\node(HH) at (-1,-4) {$abab$};
\node(EE) at (2,-3) {$bab$};
\node(EEE) at (1,-4) {$baba$};
\node(GG) at (6,-3) {$b^2a$};
\node(GGG) at (7,-4) {$b^2ab$};
\node (ZA) at (0,-5) {$\square$};
\node (ZFFF) at (-7,-5) {$a^2ba\square$};
\node (ZFF) at (-6,-5) {$a^2b\square$};
\node (ZF) at (-5,-5) {$a^2\square$};
\node (ZB) at (-4,-5) {$a\square$};
\node (ZD) at (-3,-5) {$ab\square$};
\node (ZH) at (-2,-5) {$aba\square$};
\node (ZHH) at (-1,-5) {$abab\square$};
\node (ZEEE) at (1,-5) {$baba\square$};
\node (ZEE) at (2,-5) {$bab\square$};
\node (ZE) at (3,-5) {$ba\square$};
\node (ZC) at (4,-5) {$b\square$};
\node (ZG) at (5,-5) {$b^2\square$};
\node (ZGG) at (6,-5) {$b^2a\square$};
\node (ZGGG) at (7,-5) {$b^2ab\square$};
\draw[edge,blue,thick] (A) -- (ZA) node[midway, left] {$\square$};
\draw[edge,blue,thick] (HH) -- (ZHH) node[midway, left] {$\square$};
\draw[edge,blue,thick] (H) -- (ZH) node[midway, left] {$\square$};
\draw[edge,blue,thick] (D) -- (ZD) node[midway, left] {$\square$};
\draw[edge,blue,thick] (B) -- (ZB) node[midway, left] {$\square$};
\draw[edge,blue,thick] (F) -- (ZF) node[midway, left] {$\square$};
\draw[edge,blue,thick] (FF) -- (ZFF) node[midway, left] {$\square$};
\draw[edge,blue,thick] (FFF) -- (ZFFF) node[midway, left] {$\square$};
\draw[edge,blue,thick] (EEE) -- (ZEEE) node[midway, right] {$\square$};
\draw[edge,blue,thick] (EE) -- (ZEE) node[midway, right] {$\square$};
\draw[edge,blue,thick] (E) -- (ZE) node[midway, right] {$\square$};
\draw[edge,blue,thick] (C) -- (ZC) node[midway, right] {$\square$};
\draw[edge,blue,thick] (G) -- (ZG) node[midway, right] {$\square$};
\draw[edge,blue,thick] (GG) -- (ZGG) node[midway, right] {$\square$};
\draw[edge,blue,thick] (GGG) -- (ZGGG) node[midway, right] {$\square$};
\draw[edge,blue,thick] (A) -- (B) node[midway, above] {$a$\;};
\draw[edge,thick,blue] (A) -- (C) node[midway, above] {\;$b$};
\draw[edge,thick] (B) -- (F) node[midway,right] {$a$};
\path (F) edge[->,thick, red,dashed, bend left=30] node[midway,left] {$a$} (B);
\draw[edge,thick] (F) -- (FF) node[midway,right] {$b$};
\path (FF) edge[->,thick, red, dashed, bend left=30] node[midway,left] {$b$} (F);
\draw[edge,thick] (FF) -- (FFF) node[midway,right] {$a$};
\path (FFF) edge[->,thick, red, dashed, bend left=30] node[midway,left] {$a$} (FF);
\path (FFF) edge[->,thick, red,dashed, bend left=90] node[midway,left] {$b$} (B);
\draw[edge,thick] (B) -- (D) node[midway,left] {$b$};
\path (D) edge[->,thick, red, dashed, bend right=30] node[midway,right] {$b$} (B);
\draw[edge,thick] (D) -- (H) node[midway,left] {$a$};
\path (H) edge[->,thick, red, dashed, bend right=30] node[midway,right] {$a$} (D);
\draw[edge,thick] (H) -- (HH) node[midway,left] {$b$};
\path (HH) edge[->,thick, red, dashed, bend right=30] node[midway,right] {$b$} (H);
\path (HH) edge[->,thick, red,dashed, bend right=90] node[midway,right] {$a$} (B);
\draw[edge,thick] (C) -- (E) node[midway,right] {$a$};
\path (E) edge[->,thick, red,dashed, bend left=30] node[midway,left] {$a$} (C);
\draw[edge,thick] (E) -- (EE) node[midway,right] {$b$};
\path (EE) edge[->,thick, red,dashed, bend left=30] node[midway,left] {$b$} (E);
\draw[edge,thick] (EE) -- (EEE) node[midway,right] {$a$};
\path (EEE) edge[->,thick, red,dashed, bend left=30] node[midway,left] {$a$} (EE);
\draw[edge,thick] (C) -- (G) node[midway,left] {$b$};
\path (G) edge[->,thick, red,dashed, bend right=30] node[midway,right] {$b$} (C);
\draw[edge,thick] (G) -- (GG) node[midway,left] {$a$};
\path (GG) edge[->,thick, red,dashed, bend right=30] node[midway,right] {$a$} (G);
\draw[edge,thick] (GG) -- (GGG) node[midway,left] {$b$};
\path (GGG) edge[->,thick, red,dashed, bend right=30] node[midway,right] {$b$} (GG);
\path (EEE) edge[->,thick, red,dashed, bend left=90] node[midway,left] {$b$} (C);
\path (GGG) edge[->,thick, red,dashed, bend right=90] node[midway,right] {$a$} (C);
\end{tikzpicture}
}
\end{center}
\caption{\label{figure.mccammond} The McCammond expansion of $\mathsf{KR}(S,A)$ 
of Example~\ref{example.KR}. Transition edges are blue. The edges $(p,a,q) \in E$ with $\ell(q)=\ell(p)+1$ are solid, 
whereas the edges with $\ell(q) \leqslant \ell(p)$ are dashed and red. The spanning tree $\mathsf{T}$ is obtained by
removing all the dashed red arrows.
}
\end{figure}
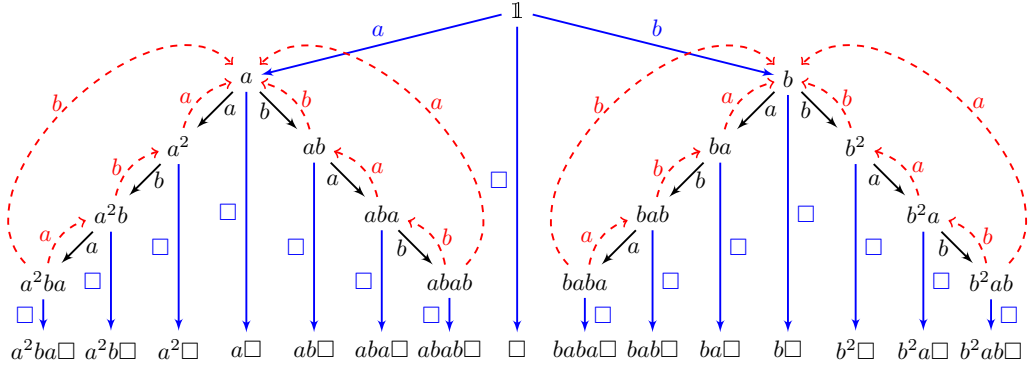

\begin{example}
\label{example.mccammond}
The McCammond expansion of $\mathsf{KR}(S,A)$ of Example~\ref{example.KR} is given in
Figure~\ref{figure.mccammond}.
\end{example}

By Remark~\ref{remark.spanning tree}, the McCammond expansion $\mathsf{Mc} \circ \mathsf{KR}(S,A)$
has a spanning tree $\mathsf{T}$. In this tree, the vertices are naturally labeled by the sequence of
edge labels in the path from $\mathbbm{1}$ to the vertex. More concretely, if
\[
	p=\left( \mathbbm{1} \stackrel{a_1}{\longrightarrow} v_1 \stackrel{a_2}{\longrightarrow} \cdots 
	\stackrel{a_\ell}{\longrightarrow} v_\ell \right)
\]
is a path in $\mathsf{T}$, then the vertex $v_\ell$ is naturally labeled by $a_1\ldots a_\ell$. Hence the corresponding
vertex $v_\ell$ in $\mathsf{Mc} \circ \mathsf{KR}(S,A)$ has a \defn{normal form} given by $a_1\ldots a_\ell$.

Remark~\ref{remark.spanning tree} also ensures that $\mathsf{Mc} \circ \mathsf{KR}(S,A)$ has the unique simple
path property, defined as follows.

\begin{definition}[Unique simple path property]
A rooted graph $(\Gamma,\mathbbm{1})$ with root $\mathbbm{1}$ has the \defn{unique simple path property}
if for each vertex $v$ in $\Gamma$ there is a unique simple path from the root  $\mathbbm{1}$ to $v$.
\end{definition}

Elements in the semaphore code $\mathcal{S}(S,A)$ are paths in $\mathsf{Mc} \circ \mathsf{KR}(S,A)$
(rather than in $\mathsf{T}$) starting at $\mathbbm{1}$ and ending in $K(S)$. They are also in natural correspondence with 
words $a_1 \ldots a_\ell \in A^+$ such that $[a_1 \cdots a_\ell]_S \in K(S)$ and $[a_1 \cdots a_{\ell-1}]_S \not \in K(S)$.
From the semaphore code, one can obtain the normal form by stripping away all loops in the path.

\subsection{Definition of $\mathsf{Pict}$}
\label{section.pict}

We are now going to define the map $\mathsf{Pict}$ from the set of tuples $(\Gamma,p)$, where $\Gamma$ is a graph
with the unique simple path property and $p$ is a simple path in $\Gamma$ starting at $\mathbbm{1}$, to the set of loop 
graphs. The straight line, that the loop graph is based on, will correspond to $p$. 
The map $\mathsf{Pict}$ was first defined by McCammond (we give a simplified definition here).

\begin{definition}[McCammond]
\label{definition.pict}
Let $\Gamma$ be a graph with the unique simple path property and $p$ a simple path in $\Gamma$ starting at
$\mathbbm{1}$. Then $\mathsf{Pict}(\Gamma,p)$ is defined by the principle of induction.

\noindent
\textbf{Induction basis:} Set $P=p$ and start at vertex $v_0=\mathbbm{1}$.

\noindent
\textbf{Induction step:} Suppose one is at vertex $v_0 \neq \tau(p)$ on path $p$. Take the edge $e$ from $v_0$ to $v_1$ 
in $p$.
\begin{enumerate}
\item If there is no edge in $\Gamma$ coming into $v_1$ besides $e$, continue with the unique next vertex in $p$,
now denoted $v_1$ (with the current vertex $v_1$ relabeled $v_0$), unless
$v_1=\tau(p)$. If $v_1=\tau(p)$, then output $\mathsf{Pict}(\Gamma,p)=P$.
\item Otherwise there is at least one edge $e' \neq e$ in $\Gamma$ going into $v_1$, given by
$e'=\left(v' \stackrel{a}{\longrightarrow} v_1\right)$ for some $a\in A$. Since $\Gamma$ has the unique simple path property by 
assumption, there must be a unique simple path starting at $\mathbbm{1}$ going to $v_0$ along the path $p$ followed by 
the path $p'$ starting at $v_0$, going along $e$ to $v_1$, and ending at $v'$.
\begin{enumerate}
\item Run the induction on $p'$ in a subgraph $\Gamma'$ of $\Gamma$, 
consisting of all edges and vertices on circuits containing a vertex of $p'$.
Note that $p'$ is simple in $\Gamma'$. The output is $P'=\mathsf{Pict}(\Gamma',p')$.
\item Modify $P$ by attaching $P'$ disjointly except at $v_1$ and adding edge $e'$ from $v'$ in $P'$ back to~$v_1$.
\end{enumerate}
\item Repeat step (2) for each edge $e'\neq e$ at vertex $v_1$.
\item Continue with the induction step unless $v_1=\tau(p)$. If $v_1=\tau(p)$, then output $\mathsf{Pict}(\Gamma,p)=P$.
\end{enumerate}
\end{definition}

\begin{remark}
If $\Gamma$ is a rooted graph with the unique simple path property, then $\Gamma$ with some edges removed (and any 
vertices that are no longer connected to the root $\mathbbm{1}$) still has the unique simple path property. This is the case 
since either the unique simple path from $\mathbbm{1}$ to $v$ is still there or the vertex $v$ is now disconnected from 
$\mathbbm{1}$ and has hence been removed.

The graph $\Gamma'$ in the Induction step (2)(a) in the definition of $\mathsf{Pict}$ can be obtained in two steps.
First remove all incoming and outgoing edges on the vertices along the path $p$ from $\mathbbm{1}$ to $v_1$, except 
the edges on the path $p$ itself. Remove all vertices that have become disconnected in this process. By the remark
above, the resulting graph still has the unique simple path property. In this graph, all simple paths go through the vertex
$v_1$. Hence we may make $v_0$ the root (removing all vertices $\mathbbm{1}$ up to $v_0$ along $p$). The result
is $\Gamma'$, which still has the unique simple path property.
\end{remark}

\begin{example}
\label{example.pict}
Let $p= \left(\mathbbm{1} \stackrel{a}{\longrightarrow} 1 \stackrel{b}{\longrightarrow} 2 \stackrel{c}{\longrightarrow} 3 \right)$ 
in
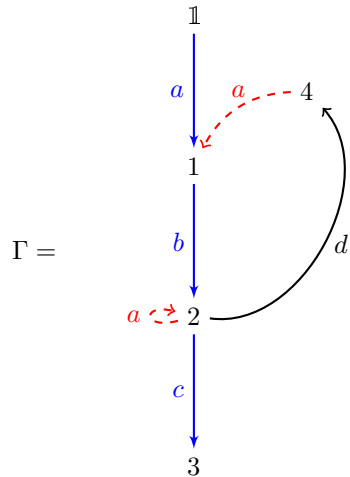
\begin{center}
\raisebox{3cm}{$\Gamma =$} \qquad  
\begin{tikzpicture}[auto]
\node (A) at (0, 0) {$\mathbbm{1}$};
\node (B) at (0,-2) {$1$};
\node(C) at (0,-4) {$2$};
\node(D) at (0,-6) {$3$};
\node (V) at (1.5,-1) {$4$};

\draw[edge,blue,thick] (A) -- (B) node[midway, left] {$a$};
\draw[edge,blue,thick] (B) -- (C) node[midway, left] {$b$};
\draw[edge,blue,thick] (C) -- (D) node[midway, left] {$c$};
\path (V) edge[->,thick, red, dashed, bend right = 30] node[midway,above] {$a$} (B);
\path (C) edge[->,thick, bend right = 70] node[midway,right] {$d$} (V);
\path (C) edge [thick,red, dashed, loop left] node {$a$} (C);
\end{tikzpicture}
\end{center}
To compute $\mathsf{Pict}(\Gamma,p)$, we start with $P=p$, $v_0=\mathbbm{1}$ and $v_1=1$.  We are in
step (2) of the Induction step with $e=\left(\mathbbm{1} \stackrel{a}{\longrightarrow} 1 \right)$ and 
$e'=\left(4\stackrel{a}{\longrightarrow} 1 \right)$. Then 
$p'=\left(\mathbbm{1} \stackrel{a}{\longrightarrow} 1 \stackrel{b}{\longrightarrow} 2 \stackrel{d}{\longrightarrow} 4 \right)$
and $\Gamma'$ is $\Gamma$ with the arrow labelled $a$ from $v'=4$ to $v_1=1$ removed. Also 
$P'=\mathsf{Pict}(\Gamma',p')$ is $p'$ with a loop labelled $a$ at vertex $2$. Attaching $P'$ at $v_1=1$ (with its vertex $2$
relabelled to $2'$ to avoid repetition) and adding edge $e'$ we obtain
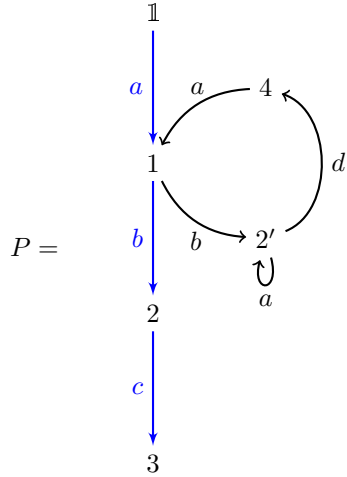
\begin{center}
\raisebox{3cm}{$P =$} \qquad  
\begin{tikzpicture}[auto]
\node (A) at (0, 0) {$\mathbbm{1}$};
\node (B) at (0,-2) {$1$};
\node(C) at (0,-4) {$2$};
\node(D) at (0,-6) {$3$};
\node (V) at (1.5,-1) {$4$};
\node (W) at (1.5,-3) {$2'$};

\draw[edge,blue,thick] (A) -- (B) node[midway, left] {$a$};
\draw[edge,blue,thick] (B) -- (C) node[midway, left] {$b$};
\draw[edge,blue,thick] (C) -- (D) node[midway, left] {$c$};
\path (V) edge[->,thick, bend right = 30] node[midway,above] {$a$} (B);
\path (B) edge[->,thick, bend right = 30] node[midway,below] {$b$} (W);
\path (W) edge[->,thick, bend right = 70] node[midway,right] {$d$} (V);
\path (W) edge [thick, loop below] node {$a$} (W);
\end{tikzpicture}
\end{center}
Since there are no further edges going into vertex $v_1=1$, we continue with the induction along $p$. This means
that we set $v_0=1$, $v_1=2$, and $e=\left(1 \stackrel{b}{\longrightarrow} 2 \right)$. Besides $e$, there is only one
other arrow going into $v_1=2$ in $\Gamma$, namely $e'=\left(2 \stackrel{a}{\longrightarrow} 2 \right)$. In this
case $p'=1 \stackrel{b}{\longrightarrow} 2$ and $\Gamma'$ is $\Gamma$ with $\mathbbm{1}$ and the arrows
$\mathbbm{1} \stackrel{a}{\longrightarrow} 1$, $4 \stackrel{a}{\longrightarrow} 1$, and $2 \stackrel{a}{\longrightarrow} 2$
removed. Hence the new $P$ with $P'=\mathsf{Pict}(\Gamma',p')$ added is
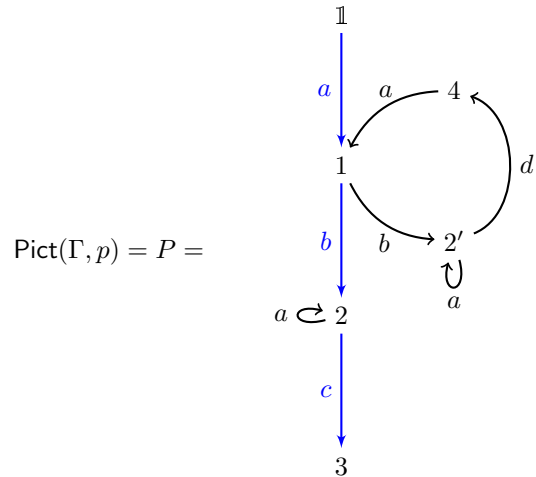
\begin{center}
\raisebox{3cm}{$\mathsf{Pict}(\Gamma,p)=P =$} \qquad  
\begin{tikzpicture}[auto]
\node (A) at (0, 0) {$\mathbbm{1}$};
\node (B) at (0,-2) {$1$};
\node(C) at (0,-4) {$2$};
\node(D) at (0,-6) {$3$};
\node (V) at (1.5,-1) {$4$};
\node (W) at (1.5,-3) {$2'$};

\draw[edge,blue,thick] (A) -- (B) node[midway, left] {$a$};
\draw[edge,blue,thick] (B) -- (C) node[midway, left] {$b$};
\draw[edge,blue,thick] (C) -- (D) node[midway, left] {$c$};
\path (V) edge[->,thick, bend right = 30] node[midway,above] {$a$} (B);
\path (B) edge[->,thick, bend right = 30] node[midway,below] {$b$} (W);
\path (W) edge[->,thick, bend right = 70] node[midway,right] {$d$} (V);
\path (W) edge [thick, loop below] node {$a$} (W);
\path (C) edge [thick, loop left] node {$a$} (C);
\end{tikzpicture}
\end{center}
The remaining induction steps do not change this $P$, which is hence also $\mathsf{Pict}(\Gamma,p)$.
\end{example}

\begin{example}
\label{example.pict1}
Consider the McCammond expansion $\Gamma= \mathsf{Mc} \circ \mathsf{KR}(S,A)$ of 
Example~\ref{example.mccammond} (see also Figure~\ref{figure.mccammond}) and the path in the 
McCammond tree $\mathsf{T}$ given by $ab\square$. Then $\mathsf{Pict}(\Gamma,ab\square)$ is
given by
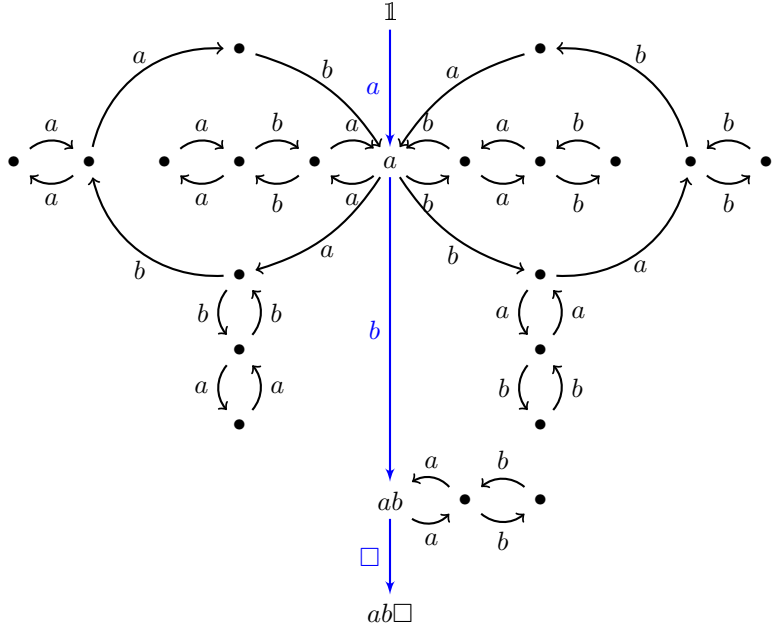
\begin{center}
\begin{tikzpicture}[auto]
\node (A) at (0, 0) {$\mathbbm{1}$};
\node (B) at (0, -2) {$a$};
\node (C) at (0, -6.5) {$ab$};
\node (D) at (0, -8) {$ab\square$};
\draw[edge,blue,thick] (A) -- (B) node[midway, left] {$a$};
\draw[edge,blue,thick] (B) -- (C) node[midway, left] {$b$};
\draw[edge,blue,thick] (C) -- (D) node[midway, left] {$\square$};
\node (P1) at (1,-6.5) {$\bullet$};
\node (P2) at (2,-6.5) {$\bullet$};
\path (C) edge[->,thick,bend right=40] node[midway,below] {$a$} (P1);
\path (P1) edge[->,thick,bend right=40] node[midway,above] {$a$} (C);
\path (P1) edge[->,thick,bend right=40] node[midway,below] {$b$} (P2);
\path (P2) edge[->,thick,bend right=40] node[midway,above] {$b$} (P1);
\node (L1) at (-1,-2) {$\bullet$};
\node (L2) at (-2,-2) {$\bullet$};
\node (L3) at (-3,-2) {$\bullet$};
\path (B) edge[->,thick,bend left=40] node[midway,below] {$a$} (L1);
\path (L1) edge[->,thick,bend left=40] node[midway,above] {$a$} (B);
\path (L1) edge[->,thick,bend left=40] node[midway,below] {$b$} (L2);
\path (L2) edge[->,thick,bend left=40] node[midway,above] {$b$} (L1);
\path (L2) edge[->,thick,bend left=40] node[midway,below] {$a$} (L3);
\path (L3) edge[->,thick,bend left=40] node[midway,above] {$a$} (L2);
\node (R1) at (1,-2) {$\bullet$};
\node (R2) at (2,-2) {$\bullet$};
\node (R3) at (3,-2) {$\bullet$};
\path (B) edge[->,thick,bend right=40] node[midway,below] {$b$} (R1);
\path (R1) edge[->,thick,bend right=40] node[midway,above] {$b$} (B);
\path (R1) edge[->,thick,bend right=40] node[midway,below] {$a$} (R2);
\path (R2) edge[->,thick,bend right=40] node[midway,above] {$a$} (R1);
\path (R2) edge[->,thick,bend right=40] node[midway,below] {$b$} (R3);
\path (R3) edge[->,thick,bend right=40] node[midway,above] {$b$} (R2);
\node (BL1) at (-2,-3.5) {$\bullet$};
\node (BL2) at (-4,-2) {$\bullet$};
\node (BL3) at (-2,-0.5) {$\bullet$};
\path (B) edge[->,thick,bend left=20] node[midway,below] {$a$} (BL1);
\path (BL1) edge[->,thick,bend left=40] node[midway,below] {$b$} (BL2);
\path (BL2) edge[->,thick,bend left=40] node[midway,above] {$a$} (BL3);
\path (BL3) edge[->,thick,bend left=20] node[midway,above] {$b$} (B);
\node (BR1) at (2,-3.5) {$\bullet$};
\node (BR2) at (4,-2) {$\bullet$};
\node (BR3) at (2,-0.5) {$\bullet$};
\path (B) edge[->,thick,bend right=20] node[midway,below] {$b$} (BR1);
\path (BR1) edge[->,thick,bend right=40] node[midway,below] {$a$} (BR2);
\path (BR2) edge[->,thick,bend right=40] node[midway,above] {$b$} (BR3);
\path (BR3) edge[->,thick,bend right=20] node[midway,above] {$a$} (B);
\node (AL) at (-5,-2) {$\bullet$};
\path (BL2) edge[->,thick,bend left=40] node[midway,below] {$a$} (AL);
\path (AL) edge[->,thick,bend left=40] node[midway,above] {$a$} (BL2);
\node (AR) at (5,-2) {$\bullet$};
\path (BR2) edge[->,thick,bend right=40] node[midway,below] {$b$} (AR);
\path (AR) edge[->,thick,bend right=40] node[midway,above] {$b$} (BR2);
\node (U1) at (-2,-4.5) {$\bullet$};
\path (BL1) edge[->,thick,bend right=40] node[midway,left] {$b$} (U1);
\path (U1) edge[->,thick,bend right=40] node[midway,right] {$b$} (BL1);
\node (U2) at (-2,-5.5) {$\bullet$};
\path (U1) edge[->,thick,bend right=40] node[midway,left] {$a$} (U2);
\path (U2) edge[->,thick,bend right=40] node[midway,right] {$a$} (U1);
\node (U3) at (2,-4.5) {$\bullet$};
\path (BR1) edge[->,thick,bend right=40] node[midway,left] {$a$} (U3);
\path (U3) edge[->,thick,bend right=40] node[midway,right] {$a$} (BR1);
\node (U4) at (2,-5.5) {$\bullet$};
\path (U3) edge[->,thick,bend right=40] node[midway,left] {$b$} (U4);
\path (U4) edge[->,thick,bend right=40] node[midway,right] {$b$} (U3);
\end{tikzpicture}
\end{center}
Following the algorithm explained in Section~\ref{section.kleene}, a Kleene expression for
$\mathcal{P}_{\mathsf{Pict}(\Gamma,ab\square)}$ is given by
\[
	L = a \{\ell_1,\ell_2,\ell_3,\ell_4\}^\star b \ell_5^\star \square,
\]
where
\begin{equation*}
\begin{split}
	\ell_1 &= a (b(aa)^\star b)^\star b (aa)^\star ab,\\
	\ell_2 &= a(b(aa)^\star b)^\star a,\\
	\ell_3 &= b (a(bb)^\star a)^\star a (bb)^\star ba,\\
	\ell_4 &= b(a(bb)^\star a)^\star b,\\
	\ell_5 &= a(bb)^\star a.
\end{split}
\end{equation*}
Hence
\begin{equation*}
\begin{split}
	\Psi_{\mathsf{Pict}(\Gamma,ab\square)} &= \frac{x_a x_b x_\square}
	{\left(1-\frac{x_a^2x_b^2}{\left(1-\frac{x_b^2}{1-x_a^2}\right)(1-x_a^2)} -\frac{x_a^2}{1-\frac{x_b^2}{1-x_a^2}}
	-\frac{x_a^2x_b^2}{\left(1-\frac{x_a^2}{1-x_b^2}\right)(1-x_b^2)} -\frac{x_b^2}{1-\frac{x_a^2}{1-x_b^2}}\right)
	\left(1-\frac{x_a^2}{1-x_b^2}\right)}\\
	&=\frac{x_a x_b x_\square(1-x_b^2)}
         {\left(1-\frac{2x_a^2x_b^2}{1-x_a^2-x_b^2} -\frac{x_a^2(1-x_a^2)}{1-x_a^2-x_b^2}
	-\frac{x_b^2(1-x_b^2)}{1-x_a^2-x_b^2}\right) (1-x_a^2-x_b^2)}\\
	&=\frac{x_a x_b x_\square(1-x_b^2)}
         {1-2x_a^2-2x_b^2+(x_a^2-x_b^2)^2}.
\end{split}
\end{equation*}
Using that $x_a+x_b+x_\square=1$, we find that in the limit $x_\square \to 0$
\[
	\lim_{x_\square\to 0} \Psi_{\mathsf{Pict}(\Gamma,ab\square)} = \frac{1}{8}(1-x_b^2).
\]
In a similar fashion, we find
\begin{equation*}
\begin{split}
	\Psi_\square & = x_\square \qquad \qquad \qquad \qquad \qquad \qquad \;\;
	 \stackrel{x_\square\to 0}{\longrightarrow} \qquad 0\\
	\Psi_{a\square} &= \frac{x_a(1-x_a^2-x_b^2)x_\square}{1-2x_a^2-2x_b^2+(x_a^2-x_b^2)^2}
	\qquad \stackrel{x_\square\to 0}{\longrightarrow} \qquad \frac{x_a}{4}\\
	\Psi_{aba\square} &= \frac{x_a^2 x_b x_\square}{1-2x_a^2-2x_b^2+(x_a^2-x_b^2)^2}
	\qquad \stackrel{x_\square\to 0}{\longrightarrow} \qquad \frac{x_a}{8}\\
	\Psi_{abab\square} &= \frac{x_a^2 x^2_b x_\square}{1-2x_a^2-2x_b^2+(x_a^2-x_b^2)^2}
	\qquad \stackrel{x_\square\to 0}{\longrightarrow} \qquad \frac{x_a x_b}{8}\\
	\Psi_{a^2\square} &= \frac{x_a^2 (1-x_a^2) x_\square}{1-2x_a^2-2x_b^2+(x_a^2-x_b^2)^2}
	\qquad \stackrel{x_\square\to 0}{\longrightarrow} \qquad \frac{x_a (1+x_a)}{8}\\
	\Psi_{a^2b\square} &= \frac{x_a^2 x_b x_\square}{1-2x_a^2-2x_b^2+(x_a^2-x_b^2)^2}
	\qquad \stackrel{x_\square\to 0}{\longrightarrow} \qquad \frac{x_a}{8}\\
	\Psi_{a^2ba\square} &= \frac{x_a^3 x_b x_\square}{1-2x_a^2-2x_b^2+(x_a^2-x_b^2)^2}
	\qquad \stackrel{x_\square\to 0}{\longrightarrow} \qquad \frac{x_a^2}{8}.
\end{split}
\end{equation*}
The stationary probabilities for the elements with $a$ and $b$ interchanged are obtained by symmetry.
It is not hard to check that these probabilities sum to one as desired.
\end{example}

As noted in the introduction, $\mathsf{Pict}(\Gamma,p)$ is not necessarily deterministic. There can be
several arrows leaving a vertex labeled by the same element $a\in A$. For example, vertex $1$ in
Example~\ref{example.pict} has two arrows labeled $b$ coming out.

One can make a non-deterministic automata $\mathcal{A}$ deterministic as follows. If $\mathcal{A}$ 
has states $Q$ with start state $\mathbbm{1}$ and final states $F$ not containing $\mathbbm{1}$, 
we make a deterministic automata $\mathsf{det}(\mathcal{A})$ accepting the same strings
going from $\mathbbm{1}$ to a member of $F$ as follows.  The states $Q'$ of $\mathsf{det}(\mathcal{A})$ are the 
collection of subsets of $Q$ determined a follows:  
\begin{itemize}
\item $\{\mathbbm{1}\}$ is in $Q'$;
\item if $Z\in Q'$, then $Z.a\in Q'$ for $a\in A$, where $Z.a= \{q \mid z \stackrel{a}{\longrightarrow} q \in \mathcal{A} 
\text{ where } z\in Z\}$.
\end{itemize}
One continues by induction until the process adds no new subsets. For $\mathsf{det}(\mathcal{A})$, start in state 
$\{\mathbbm{1}\}$. The final states are all the states of $\mathsf{det}(\mathcal{A})$ such that the intersection with $F$
is non-empty.

With this definition, making $\mathsf{Pict}(\Gamma,p)$ deterministic gives the automata for $(\Gamma,p)$ back.

\subsection{Proof of Theorem~\ref{theorem.main}}
\label{section.proof}

As explained in Section~\ref{section.markov}, any finite Markov chain $\mathcal{M}$ can be described
as a Markov chain $\mathcal{M}(S,A)$ in terms of a finite semigroup $S$ with generators $A$. Since by
Theorem~\ref{theorem.stationary real}, $\Psi_w^{{\mathcal{M}(S,A)}}$ is the sum over 
$\Psi_v^{{\mathcal{M}(\mathsf{KR}(S,A))}}$, it suffices to prove the statement of Theorem~\ref{theorem.main}
for $\Psi_v^{{\mathcal{M}(\mathsf{KR}(S,A))}}$. When $K(S)$ is not left zero,
we may use the limiting construction of~\eqref{equation.limit} to obtain $\Psi_v^{{\mathcal{M}(\mathsf{KR}(S,A))}}$
from the case in which the minimal ideal is left zero. Assuming that $K(S)$ is left zero, we have
by Theorem~\ref{theorem.stationary general}
\begin{equation}
\label{equation.psi semaphore}
	\Psi^{\mathcal{M}(\mathsf{KR}(S,A))}_w 
	= \sum_{\substack{s \in \mathcal{S}(S,A)\\ [s]_{\mathsf{KR}(S,A)}=w}}\; \prod_{a\in s} x_a
	\qquad \text{for all $w \in K(\mathsf{KR}(S,A))$.}
\end{equation}

As explained in Section~\ref{section.semaphore}, there is a normal form associated to each
semaphore code element $s\in \mathcal{S}(S,A)$. Namely, $s$ is a path in $\mathsf{Mc} \circ \mathsf{KR}(S,A)$
starting at $\mathbbm{1}$ and the normal form is the simple path with all loops stripped away from $s$; equivalently 
the normal form is the path in $\mathsf{T}$ starting at $\mathbbm{1}$ and ending at $\tau(s)$, where $\mathsf{T}$
is the tree associated to the McCammond expansion $\mathsf{Mc} \circ \mathsf{KR}(S,A)$.
In the tree $\mathsf{T}$, a path $p$ starting at $\mathbbm{1}$ is also naturally in bijection with its endpoint $\tau(p)$. 
Hence we may identify vertex $t\in \mathsf{T}$ with the path from $\mathbbm{1}$ to $t$ in $\mathsf{T}$ or equivalently with 
the simple path from $\mathbbm{1}$ to $t$ in $\mathsf{Mc} \circ \mathsf{KR}(S,A)$. Therefore, we may rewrite the sum 
in~\eqref{equation.psi semaphore} as
\begin{equation}
\label{equation.psi T}
	\Psi^{\mathcal{M}(\mathsf{KR}(S,A))}_w 
	= \sum_{\substack{t \in \mathsf{T}\\ [t]_{\mathsf{KR}(S,A)}=w}}\; 
	\left(\sum_{\substack{s\in \mathcal{S}(S,A)\\ \tau(s)=t}} \; \prod_{a\in s} x_a\right)
	\qquad \text{for all $w \in K(\mathsf{KR}(S,A))$.}
\end{equation}
We claim that for a given $t\in \mathsf{T}$ with $[t]_{\mathsf{KR}(S,A)} \in K(\mathsf{KR}(S,A))$
\begin{equation}
\label{equation.psi pict}
	\Psi_{\mathsf{Pict}(\mathsf{Mc} \circ \mathsf{KR}(S,A),t)} 
	= \sum_{\substack{s\in \mathcal{S}(S,A)\\ \tau(s)=t}} \prod_{a\in s} x_a.
\end{equation}
Recall that by Definition~\ref{definition.normal}
\[
	\Psi_{\mathsf{Pict}(\mathsf{Mc} \circ \mathsf{KR}(S,A),t)} 
	= \sum_{p \in \mathcal{P}_{\mathsf{Pict}(\mathsf{Mc} \circ \mathsf{KR}(S,A),t)}} \prod_{a\in p} x_a.
\]
Hence~\eqref{equation.psi pict} can be proved by establishing a bijection
\begin{equation}
\label{equation.bijection}
	\varphi \colon \{s\in \mathcal{S}(S,A) \mid \tau(s)=t\} \longrightarrow
	\mathcal{P}_{\mathsf{Pict}(\mathsf{Mc} \circ \mathsf{KR}(S,A),t)}.
\end{equation}
In fact, we are going to prove a slight generalization of~\eqref{equation.bijection}. Namely, for any $t\in \mathsf{T}$
we will show that there is a bijection
\begin{equation}
\label{equation.bijection 1}
	\varphi \colon \{s\in \mathcal{P}_{\mathsf{Mc} \circ \mathsf{KR}(S,A)} \mid \tau(s)=t\} \longrightarrow
	\mathcal{P}_{\mathsf{Pict}(\mathsf{Mc} \circ \mathsf{KR}(S,A),t)},
\end{equation}
where $\mathcal{P}_{\mathsf{Mc} \circ \mathsf{KR}(S,A)}$ is the set of paths in
$\mathsf{Mc} \circ \mathsf{KR}(S,A)$ starting at $\mathbbm{1}$. Then~\eqref{equation.bijection} is the 
special case when $[t]_{\mathsf{KR}(S,A)} \in K(\mathsf{KR}(S,A))$.

To define $\varphi$ in~\eqref{equation.bijection 1}, fix $t=a_1 \cdots a_k$, where $a_i\in A$ are the labels in the path in 
$\mathsf{T}$. A path $s \in \mathcal{P}_{\mathsf{Mc} \circ \mathsf{KR}(S,A)}$ with $\tau(s) = t$, can be viewed as 
$t$ with circuits $\ell_j^{(j)}$ interspersed. More precisely,
\[
	s=a_1\left(\prod_{j\in J_1} \ell_1^{(j)}\right) a_2 \left(\prod_{j\in J_2} \ell_2^{(j)} \right) \cdots a_k
	\left(\prod_{j\in J_k} \ell_k^{(j)}\right),
\]
where $\tau(a_1 \cdots a_i) = \tau(a_1 \cdots a_i \ell_i^{(j)})$ for all $1\leqslant i\leqslant k$ and $j\in J_i$ and any initial
subsequence of $\ell_i^{(j)}$ does not reach the vertex $a_1\cdots a_i$. Here the sets $J_i$ index the set of circuits 
$\{\ell_i^{(j)} \mid j\in J_i\}$ at vertex $a_1 \cdots a_i$ and either $J_i=\{1,2,\ldots,n_i\}$ is a finite set or $J_i=\{1,2,3,\ldots\}$
is the set of positive integers. In other words, each $\ell_i^{(j)}$ is a circuit from vertex $a_1 \cdots a_i$ to itself, which
does not pass through $a_1 \cdots a_i$ otherwise. The last step of $\ell_i^{(j)}$ is an edge in 
$\mathsf{Mc} \circ \mathsf{KR}(S,A)$ that is not in $\mathsf{T}$. 
Suppose by induction that 
\[
	s'=a_1\left(\prod_{j\in J_1} \ell_1^{(j)}\right) \cdots a_i \left(\prod_{j\in J'_i} \ell_i^{(j)}\right),
\]
where $1\leqslant i \leqslant k$ and $J'_i =\{1,2,\ldots,n_i'\} \subseteq J_i$ or $J_i'=J_i$, is mapped to $\pi$ in
$\mathsf{Pict}(\mathsf{Mc} \circ \mathsf{KR}(S,A),a_1 \cdots a_i)$ under $\varphi$. We need to distinguish two cases.

\smallskip

\noindent
\textbf{Case $J_i' \subsetneq J_i$.} Let $j$ be the smallest element in $J_i \setminus J_i'$. Recall that 
$\mathsf{Mc} \circ \mathsf{KR}(S,A)$ has the unique simple path property. Hence the path $p'$ in 
$\mathsf{Mc} \circ \mathsf{KR}(S,A)$ from $v_0=a_1 \cdots a_{i-1}$
through $v_1=a_1 \cdots a_i$ to $v'$, which is $a_1 \cdots a_i \ell_i^{(j)}$ with the last edge $e'$ removed is a path in 
$\Gamma'$ in the notation of Section~\ref{section.pict}. By induction this path is mapped to $\pi'$ in 
$\mathcal{P}_{\mathsf{Pict}(\Gamma',p')}$. Hence
\[
	\varphi(s'\ell_i^{(j)})=\pi\pi' \in \mathcal{P}_{\mathsf{Pict}(\mathsf{Mc} \circ \mathsf{KR}(S,A),a_1 \cdots a_i)}
\]
This corresponds to the induction step (2) in Definition~\ref{definition.pict}.

\smallskip

\noindent
\textbf{Case $J_i' =J_i$.} If $i=k$, we are done. If $i<k$, we define
\[
	\varphi(s'a_{i+1}) = \pi a_{i+1} \in \mathcal{P}_{\mathsf{Pict}(\mathsf{Mc} \circ \mathsf{KR}(S,A),a_1 \cdots a_{i+1})},
\]
which is a well-defined path since the last step is along the straight line path and hence unique.
This corresponds to the induction step (1) (if $J_i=\emptyset$) or step (4) (if $J_i\neq \emptyset$)
in Definition~\ref{definition.pict}.

This shows that $\varphi$ is a well-defined map. It has an inverse $\varphi^{-1}$ by mapping a path 
$\pi \in \mathcal{P}_{\mathsf{Pict}(\mathsf{Mc} \circ \mathsf{KR}(S,A),t)}$ to a path in $\mathsf{Mc} \circ \mathsf{KR}(S,A)$
by just reading the labels of the edges. This indeed gives a path in $\mathsf{Mc} \circ \mathsf{KR}(S,A)$ by the construction
of $\mathsf{Pict}$. 

Combining~\eqref{equation.psi T} and~\eqref{equation.psi pict}, we obtain
\[
	\Psi^{\mathcal{M}(\mathsf{KR}(S,A))}_w 
	=  \sum_{\substack{t \in \mathsf{T}\\ [t]_{\mathsf{KR}(S,A)}=w}}\; \Psi_{\mathsf{Pict}(\mathsf{Mc} \circ \mathsf{KR}(S,A),t)},
\]
which proves Theorem~\ref{theorem.main} since $\mathsf{Pict}(\mathsf{Mc} \circ \mathsf{KR}(S,A),t)$ is a loop graph.

In summary, we proved the following theorem, which is a more detailed version of Theorem~\ref{theorem.main}.

\begin{theorem}
\label{theorem.main 1}
Let $\mathcal{M}(S,A)$ be a Markov chain associated to the finite semigroup with generators $(S,A)$.
If $K(S)$ is left zero, the stationary distribution is given by
\[
	\Psi_w^{\mathcal{M}(S,A)} = \sum_{\substack{t \in \mathsf{T}\\ [t]_S = w}} 
	\Psi_{\mathsf{Pict}(\mathsf{Mc} \circ \mathsf{KR}(S,A),t)} \qquad \text{for $w\in K(S)$,}
\]
where $\mathsf{T}$ is the spanning tree of $\mathsf{Mc} \circ \mathsf{KR}(S,A)$.
Otherwise
\[
	\Psi_w^{\mathcal{M}(S,A)} = \sum_{\substack{t \in \mathsf{T}\\ [t]_S = w \square}} 
	\lim_{x_\square \to 0}  \Psi_{\mathsf{Pict}(\mathsf{Mc} \circ \mathsf{KR}(S \cup \{\square\},A \cup \{\square\}),t)}
	\qquad \text{for $w\in K(S)$,}
\]
where $\mathsf{T}$ is the spanning tree of $\mathsf{Mc} \circ \mathsf{KR}(S \cup \{\square\},A \cup \{\square\})$
and $\square$ acts as zero.
\end{theorem}


\bibliographystyle{alpha}
\bibliography{loop}{}

\end{document}